\documentclass[12pt,reqno,twoside]{amsart}%
\usepackage{amsmath,amsthm,amscd,amsthm,upref,indentfirst}
\usepackage{amsfonts,mathrsfs}
\usepackage{amssymb,amsbsy,bm}
\usepackage{graphicx}
\usepackage{amsaddr}
\usepackage{soul}

\usepackage{mathabx}
\usepackage[usenames,dvipsnames]{color}

\theoremstyle{plain}
\newtheorem{theorem}{Theorem}

\theoremstyle{definition}

\theoremstyle{remark}

\numberwithin{equation}{section}
\numberwithin{theorem}{section}

\normalfont\upshape
\usepackage{exscale}
\usepackage[scaled=0.86]{helvet}

\usepackage[scr=euler,scrscaled=1.1,bb=txof,bbscaled=1.16,cal=boondox,calscaled=1.1]{mathalfa}

\renewcommand{\mathbf}{\bm}
\renewcommand{\mathit}[1]{\mathscr #1}

\renewcommand{\mathfrak}[1]{{\textbf{\upshape #1}}}

\renewcommand{\mathtt}[1]{\scalebox{1.2}{\bf \texttt{\upshape#1}}}
\renewcommand{\emph}[1]{\textcolor{blue}{\textbf{#1}}}
\renewcommand{\mathrm}[1]{\scalebox{1.15}{\textbf{\upshape #1}}}
\usepackage{float}

\usepackage[justification=centering]{caption}

\numberwithin{equation}{section}
\numberwithin{theorem}{section}

\usepackage[normalem]{ulem}
\usepackage[square,comma]{natbib}

\usepackage{mparhack}

\usepackage{keyval}
\usepackage{totcount}
\newtotcounter{citnum} %From the package documentation
\def\oldbibitem{} \let\oldbibitem=\bibitem
\def\bibitem{\stepcounter{citnum}\oldbibitem}

\usepackage[verbose]{backref}
\backrefsetup{verbose=false}
\renewcommand*{\backref}[1]{}
\renewcommand*{\backrefalt}[4]{[{\tiny%
    \ifcase #1 \textsl{Not cited}%
          \or \textsl{Cited on page}~\textcolor{BrickRed}{#2}%
          \else \textsl{Cited on pages}~\textcolor{BrickRed}{#2}%
    \fi%
    }]}

\usepackage{epigraph}

\usepackage{fancyhdr}
\author{\small\scshape S\lowercase{teven} D\lowercase{uplij}}

\address{% \small \scshape
Center for Information Technology,
University of M\"unster,
48149 M\"unster,
Germany
\\and\\
Yantai Research Institute,
Harbin Engineering University,
265615 Yantai, China
}
\email{\small \sf douplii@uni-muenster.de;
sduplij@gmail.com;
http://www.uni-muenster.de/IT.StepanDouplii}

\title{\large\bfseries\scshape
P\lowercase{ositional numeral systems over polyadic rings}}

\date{\textit{of start} May 25, 2025. \textit{Date}:
\textit{of completion}
June 13, 2025.
\newline
\mbox{}\hskip 1.16em
\textit{Total}:
11
references.
}

\renewcommand{\refname}{\textsc{References}}

\let\origsection\section
\renewcommand{\section}[1]{\sectionmark{#1}\origsection{#1}}
\let\origsubsection\subsection
\renewcommand{\subsection}[1]{\subsectionmark{#1}\origsubsection{#1}}

\makeatletter
\renewenvironment{thebibliography}[1]{%
  \@xp\origsection\@xp*\@xp{\refname}%
  \normalfont\footnotesize\labelsep .9em\relax
  \renewcommand\theenumiv{\arabic{enumiv}}\let\p@enumiv\@empty
  \vspace*{-5pt}% NEW
  \list{\@biblabel{\theenumiv}}{\settowidth\labelwidth{\@biblabel{#1}}%
    \leftmargin\labelwidth \advance\leftmargin\labelsep
    \usecounter{enumiv}}%
  \sloppy \clubpenalty\@M \widowpenalty\clubpenalty
  \sfcode`\.=\@m
}{%
  \def\@noitemerr{\@latex@warning{Empty `thebibliography' environment}}%
  \endlist
}
\makeatother

\subjclass[2010]{11-03, 11A07, 11A67, 20N15}
\keywords{positional number system, arity, polyadic structure, polyadic number, congruence class, numeral, polyadic ring, mixed-base, querelement}

\begin{document}
\mbox{}
\vskip 1.8cm
\begin{abstract}
%\noindent
%TCIDATA{OutputFilter=latex2.dll}
%TCIDATA{Version=5.50.0.2953}
%TCIDATA{LaTeXparent=0,0,example.TEX}

\noindent We construct positional numeral systems that work natively over
nonderived polyadic $\left(  m,n\right)  $-rings whose addition takes $m$
arguments and multiplication takes $n$. In such rings, the length of an
admissible additive word and a multiplicative tower are not arbitrary (as in
the binary case), but \textquotedblleft quantized\textquotedblright. Our main
contributions are the following. Existence: every commutative $\left(
m,n\right)  $-ring admits a base-$p$ place-value expansion that respects the
word length constraint in terms of numbers of operation compositions
$\ell_{mult}=\ell_{add}(m-1)+1$. Lower bound: the minimum number of digits
is greater than or equal to the arity of addition $m$. Representability gap: for $m,n\geq3$ only a
proper subset of ring elements possess finite expansions, characterized by
congruence-class arity shape invariants $I^{(m)}$ and $J^{(n)}$. Mixed-base
\textquotedblleft polyadic clocks\textquotedblright: allowing a different base
at each position enlarges the design space quadratically in the digit count.
Catalogues: explicit tables for the integer rings $\mathbb{Z}_{4,3}$ and
$\mathbb{Z}_{6,5}$ illustrate how ordinary integers lift to distinct polyadic
variables. These results lay the groundwork for faster arity-aware arithmetic,
exotic coding schemes, and hardware that exploits operations beyond the binary pair.
\end{abstract}

\maketitle

\thispagestyle{empty}
\mbox{}
\vspace{-0.5cm}
%\begin{small}
\tableofcontents
%\end{small}%
\newpage

\pagestyle{fancy}

\addtolength{\footskip}{15pt}

\renewcommand{\sectionmark}[1]{%
\markboth{
%\textmd{\  \thesection.}
{ \scshape #1}}{}}

\renewcommand{\subsectionmark}[1]{%
\markright{
\mbox{\;}\\[5pt]
\textmd{#1}}{}}

\fancyhead{}
\fancyhead[EL,OR]{\leftmark}
\fancyhead[ER,OL]{\rightmark}
\fancyfoot[C]{\scshape -- \textcolor{BrickRed}{\thepage} --}
\fancyfoot[R]{}
%\fancyfoot[L]{{\small \hours}}

\renewcommand\headrulewidth{0.5pt}
\fancypagestyle {plain1}{ %
\fancyhf{}
\renewcommand {\headrulewidth }{0pt}
\renewcommand {\footrulewidth }{0pt}
}

\fancypagestyle{plain}{ %
\fancyhf{}
\fancyhead[C]{\scshape S\lowercase{teven} D\lowercase{uplij} \hskip 0.7cm \MakeUppercase{Polyadic Hopf algebras and quantum groups}}
\fancyfoot[C]{\scshape - \thepage  -}
\renewcommand {\headrulewidth }{0pt}
\renewcommand {\footrulewidth }{0pt}
}

\fancypagestyle{fancyref}{ %
\fancyhf{} % remove everything
\fancyhead[C]{\scshape R\lowercase{eferences} }
\fancyfoot[C]{\scshape -- \textcolor{BrickRed}{\thepage} --}
\renewcommand {\headrulewidth }{0.5pt}
\renewcommand {\footrulewidth }{0pt}
}

\fancypagestyle{emptyf}{
\fancyhead{}
\fancyfoot[C]{\scshape -- \textcolor{BrickRed}{\thepage} --}
\renewcommand{\headrulewidth}{0pt}
}
%\mbox{}
%\vskip 1cm
\thispagestyle{emptyf}

\section{\textsc{Introduction}}

Classical positional numeral systems (for a review see, e.g.
\cite{con/guy,ifrah}) answer the question: how to encode quantities
efficiently? The linked question, which discusses how to operate on those
encodings, has the narrow answer: the admissible operations are binary, that
is addition and multiplication decompose into chains of pairwise interactions.
Nowadays in group theory, abstract algebra, coding theory, ultra-metric
analysis and physics $n$-ary (or polyadic) operations, those that swallow $n$
inputs at once, crop up with growing frequency (for a review, see
\cite{duplij2022} and refs therein), nevertheless their foundational
connections to numeral systems have remained largely unexplored (with the
exception of \cite{dup2017a}). This paper fills that gap and bridges two
seemingly disparate domains: the positional numeral systems that represent
integers \cite{men/bro,fomin74} and polyadic rings that generalize ring theory
to $n$-ary operations \cite{lee/but}.

The main novel and unexpected phenomena which were uncovered in this paper are
the following:

Double quantization of word length. As a starting point is the rich body of
work on polyadic algebraic structures, we observe that in an $n$-ary magma the
fundamental operation $\mu_{n}$ combines $n$ elements at a time, such that
closeness and total associativity now \textquotedblleft
quantizes\textquotedblright\ the admissible word lengths to the arithmetic
progression (opposed to the binary case, where word length is arbitrary).
Adding a second $m$-operation $\nu_{m}$ yields an $\left(  m,n\right)  $-ring,
provided the two interact through a generalized distributive law and obey
\textquotedblleft double quantization\textquotedblright\ of word length. The
polyadic rings enjoy exotic features: multiple identities, two kinds of
querelements replacing inverses (of addition and multiplication), and the
possibility that zero or one may fail to exist at all.

Numeral polyadization. We introduce a novel systematic recipe that translates
the classical place-value expansion into a special genuine composition of
$m$-ary additions and $n$-ary multiplications which reduces to the standard
positional numeral presentation in the binary $\left(  2,2\right)  $ case. The
key feature of polyadic numeral presentation is the unforeseen fact that the
number of summands and the number of multiplicative factors are no longer
independent, but lock together through the arities, forcing numeral strings to
respect a built-in \textquotedblleft polyadic uncertainty
principle\textquotedblright.

Numeral representability is no longer guaranteed. Unlike the binary case
($m=n=2$), not every element of an $(m,n)$-ring admits a numeral expansion. We
arrive to the theorem claiming that for $m,n\geq3$, the set of numbers with a
finite positional expansion becomes a strict subset (characterized by the
shape invariants $I^{\left(  m\right)  }$ and $J^{\left(  n\right)  }$) of the
$\left(  m,n\right)  $-ring, unless $m=n=2$.

Mixed-base and \textquotedblleft polyadic clock\textquotedblright\ systems.
Generalizing time-of-day notation, we allow each digit position to carry its
own $n$-ary tower of bases, yielding a combinatorial explosion of admissible
encodings and opening the door to cryptographic and coding-theoretic
applications, as well as turn coin-change, currency, and time-keeping problems
to higher arity.

In general, polyadic numeral systems promise new tools wherever data are
naturally grouped in arity-specific chunks: multi-qubit gates in quantum
circuits, multiary neural activations, and even higher-order interaction terms
in field theory.

The structure of the paper is as follows. \textsc{Section} $2$ recaps the
essentials of polyadic one set algebraic structures with an emphasis on the
quantization of admissible word lengths and the various exotic element types
(querelements, polyadic powers and idempotents, etc.). \textsc{Section} $3$
develops the theory of $\left(  m,n\right)  $-rings built from integer
representatives of congruence classes, culminating in the arity-shape map that
dictates which classes admit the polyadic ring structure. \textsc{Section} $4$
our main contribution constructs positional numeral systems over the $\left(
m,n\right)  $-rings, proves the minimal-digit theorem, and works through
explicit examples in the $\left(  4,3\right)  $ and $\left(  6,5\right)  $
settings. We close with a catalogue of open problems, ranging from polyadic
floating-point arithmetic to potential links with higher-categories.

\section{\label{sec-prel}\textsc{Preliminaries}}

For self-consistency, we briefly recall polyadic structure notation and
general properties of polyadic rings (for details, see \cite{duplij2022}).

Let $X$ be a non-empty set and $X^{\times n}$ be its Cartesian product
consisting of elements $\left(  x_{1},\ldots,x_{n}\right)  $ called polyads
(or $n$-tuples) and denoted by $\left(  \mathfrak{x}\right)  \in X^{\times n}%
$, and an $n$-tuple with equal elements is $\left(  x^{n}\right)  $. On
$X^{\times n}$ the polyadic ($n$-ary operation) is defined by $\mathbf{\mu
}_{n}:X^{\times n}\longrightarrow X$ and is denoted as $\mathbf{\mu}%
_{n}\left[  \mathfrak{x}\right]  $. A polyadic structure $\left\langle
X\mid\mathbf{\mu}_{n_{i}}\right\rangle $ is the set $X$ that is closed with
respect to polyadic operations $\mathbf{\mu}_{n_{i}}$ (see, e.g.,
\cite{kurosh,cohn}).

The basic one-operation polyadic structure is called a ($n$-ary) magma
$\mathcal{M}=\left\langle X\mid\mathbf{\mu}_{n}\right\rangle $. The (totally)
commutative $n$-ary magma $\mathcal{M}=\left\langle X\mid\mathbf{\mu}_{n}\mid
comm\right\rangle $ is defined by $\mathbf{\mu}_{n}=\mathbf{\mu}_{n}%
\circ\sigma$, where $\sigma$ is an arbitrary $n$-permutation. Additional
requirements lead to various structures called group-like ones. For instance,
polyadic associative magma is an $n$-ary semigroup $\mathcal{S}_{n}%
=\left\langle X\mid\mathbf{\mu}_{n}\mid assoc\right\rangle $. The polyadic
associativity can be defined as the invariance \cite{dup2018a}%
\begin{equation}
\mathbf{\mu}_{n}\left[  \mathfrak{x},\mathbf{\mu}_{n}\left[  \mathfrak{y}%
\right]  ,\mathfrak{z}\right]  =invariant, \label{mi}%
\end{equation}
where the internal multiplication can be on any of $n-1$ places (giving $n-1$
relations) and the polyads $\mathfrak{x},\mathfrak{y},\mathfrak{z}$ are of
needed size, the total number of elements in (\ref{mi}) is fixed $2n-1$ (cf.
with the binary associativity $\left(  x\cdot y\right)  \cdot z=x\cdot\left(
y\cdot z\right)  $). This allows us to omit internal bracket in compositions
and define the iterated multiplication%
\begin{equation}
\mathbf{\mu}_{n}^{\circ\ell_{\mu}}\left[  \mathfrak{x}\right]  =\overset
{\ell_{\mu}}{\overbrace{\mathbf{\mu}_{n}\left[  \mathbf{\mu}_{n}\left[
\ldots\mathbf{\mu}_{n}\left[  \mathfrak{x}\right]  \right]  \right]  }%
},\ \ \ \mathfrak{x}\in X^{\ell_{\mu}\left(  n-1\right)  +1}, \label{ml}%
\end{equation}
where $\ell_{\mu}$ is the \textquotedblleft number\textquotedblright\ of
operations (now multiplications). From (\ref{ml}) it follows the main and
crucial peculiarity of polyadic structures which differs them from ordinary
binary ones ($n=2$): the length of a word $w_{\mu}\left(  n\right)  $ inside
composition of an $n$-ary multiplications is not arbitrary, but it is
\textquotedblleft quantized\textquotedblright, obeying the following
admissible values%
\begin{equation}
w_{\mu}^{admiss}\left(  n\right)  =\ell_{\mu}\left(  n-1\right)  +1,
\label{wl}%
\end{equation}
which means that we can multiply $w_{\mu}^{admiss}\left(  n\right)  $ elements only.

Besides, one can treat the l.h.s. of (\ref{wl}) as another operation (the
iterated multiplication) having higher arity%
\begin{equation}
\mathbf{\bar{\mu}}_{\ell_{\mu}\left(  n-1\right)  +1}=\mathbf{\mu}_{n}%
^{\circ\ell_{\mu}}\left[  \mathfrak{x}\right]  .
\end{equation}

This allows us to divide all polyadic operations into two classes: iterated
from lower arity operations and noniterated, or derived and nonderived. Such
division can also depend on the underlying set.

For simplest example, the ternary operation $\mathbf{\mu}_{3}\left[
x,y,z\right]  =xyz$ (product of three integers from $\mathbb{Z}$) is derived
for positive integers $X=\mathbb{Z}^{+}$ (because $\mathbf{\mu}_{3}\left[
x,y,z\right]  =\mathbf{\mu}_{2}\left[  \mathbf{\mu}_{2}\left[  x,y\right]
,z\right]  =\left(  x\cdot y\right)  \cdot z$), but $\mathbf{\mu}_{3}$ is
nonderived for negative integers $X=\mathbb{Z}^{-}$, because the product of
two negative integers is positive, and so the binary operation $\mathbf{\mu
}_{2}$ is not closed, i.e. the iteration does not exist. The ternary
associativity is obvious in both cases, and therefore the polyadic structures
$\mathcal{S}_{3}^{+}=\left\langle \mathbb{Z}^{+}\mid\mathbf{\mu}_{3}\mid
assoc\right\rangle $ and $\mathcal{S}_{3}^{-}=\left\langle \mathbb{Z}^{-}%
\mid\mathbf{\mu}_{3}\mid assoc\right\rangle $ are ternary semigroups,
nevertheless only the second one $\mathcal{S}_{3}^{-}$ is nonderived. Below we
will consider barely the nonderived polyadic operations and the corresponding
polyadic structures. Recall definitions of important elements in polyadic
structures which will be used later on.

For each element $x\in X$ we can construct its $\ell_{\mu}$-polyadic power by%
\begin{equation}
x^{\left\langle \ell_{\mu}\right\rangle }=\mathbf{\mu}_{n}^{\circ\ell_{\mu}%
}\left[  x^{\ell_{\mu}\left(  n-1\right)  +1}\right]  , \label{xl}%
\end{equation}
and in the binary case $n=2$ the polyadic power differs by $1$ from the
ordinary power $x^{\left\langle \ell_{\mu}\right\rangle }=x^{\ell_{\mu}+1}$.

A polyadic idempotent $x_{id}$ (if exists) is%
\begin{equation}
x_{id}^{\left\langle \ell_{\mu}\right\rangle }=x_{id},\ \ \ x_{id}\in X.
\label{id}%
\end{equation}

A polyadic zero is defined unambiguously by the $n-1$ relations%
\begin{equation}
\mathbf{\mu}_{n}\left[  z,\mathfrak{x}\right]  =z,\ \ \ \mathfrak{x}\in
X^{n-1}, \label{z}%
\end{equation}
where $z$ can be on any of $n-1$ places. A polyadic nilpotent element
$x_{nil}$ is defined by one relation%
\begin{equation}
x_{nil}^{\left\langle \ell_{\mu}\right\rangle }=z,\ \ \ x_{nil}\in X.
\label{mx}%
\end{equation}

The neutral $\left(  n-1\right)  $-polyad $\mathfrak{e}$ (in general
non-unique) is given by%
\begin{equation}
\mathbf{\mu}_{n}\left[  x,\mathfrak{e}\right]  =x,\ \ \ \mathfrak{e}\in
X^{n-1}. \label{ee}%
\end{equation}
If all elements in the neutral polyad are equal $\mathfrak{e}=e^{n-1}$, then%
\begin{equation}
\mathbf{\mu}_{n}\left[  x,e^{n-1}\right]  =x, \label{e1}%
\end{equation}
and $e$ can be on any of $n-1$ places and is called an identity of
$\left\langle X\mid\mathbf{\mu}_{n}\right\rangle $.

If the identity $e$ exists, we can define polyadic $\ell_{\mu}$-reflection by%
\begin{equation}
x_{refl}^{\left\langle \ell_{\mu}\right\rangle }=e,\ \ \ x_{refl}\in X.
\end{equation}

As follows from (\ref{z}) with $\mathfrak{x}=z^{n-1}$ and (\ref{e1}) with
$x=e$, polyadic zero $z$ and identity $e$ are idempotents satisfying
(\ref{id}). There are exotic polyadic structures without idempotents, zero,
identity at all or with several identities.

For instance, the above ternary semigroup of negative integers $\mathcal{S}%
_{3}^{-}$ contains two neutral sequences $\mathfrak{e}_{+}=\left(
+1,+1\right)  $, $\mathfrak{e}_{-}=\left(  -1,-1\right)  $ and two ternary
identities $e_{+}=+1$ and $e_{-}=-1$.

The invertibility in the polyadic case ($n\geq3$) is not connected with
identity (\ref{e1}), but is guided by the polyadic analogue of inverse, the
querelement $\bar{x}=\bar{x}\left(  x\right)  $ defined by the $n-1$ relations
\cite{dor3}%
\begin{equation}
\mathbf{\mu}_{n}\left[  \bar{x},x^{n-1}\right]  =x,\ \ \ x\in X, \label{mq}%
\end{equation}
which should be satisfied with $\bar{x}$ be on all of $n-1$ places, and such
$x$ is called polyadically invertible. The polyadic solvability is defined as
the unique solution for $x_{0}$ with given $\mathfrak{u},y$ of the $n-1$
equations%
\begin{equation}
\mathbf{\mu}_{n}\left[  x_{0},\mathfrak{u}\right]  =y,\ \ \ x_{0},y\in
X,\ \ \ \mathfrak{u}\in X^{n-1} \label{my}%
\end{equation}
with $x_{0}$ be on any of all $n-1$ places. If each element of an $n$-ary
semigroup $\mathcal{S}_{n}$ is polyadically invertible, or equivalently,
(\ref{my}) is solvable on any place, then $\mathcal{S}_{n}$ becomes an $n$-ary
(polyadic) group $\mathcal{G}_{n}=\left\langle X\mid\mathbf{\mu}_{n}\mid
assoc\mid solv\right\rangle $. The existence of identity is not necessary for
polyadic groups.

For example, the imaginary fractions $ix/y$ with $x,y\in\mathbb{Z}^{odd}$
($i^{2}=-1$, operations are in $\mathbb{C}$) form a nonderived ternary group
$\mathcal{G}_{3}^{odd}$ without identity, each element is polyadically
invertible and has its unique querelement, that is $\overline{\left(
ix/y\right)  }=-iy/x$.

\section{\label{sec-polring}\textsc{Polyadic rings}}

A set with two polyadic operations belongs to the so-called ring-like polyadic
structures (see, e.g., \cite{lee/but,duplij2022} and refs. therein).

A polyadic or $\left(  m,n\right)  $-ring $\mathcal{R}_{m,n}=\left\langle
X\mid\mathbf{\nu}_{m},\mathbf{\mu}_{n}\right\rangle $ is a nonempty set $X$
with $m$-ary addition $\mathbf{\nu}_{m}:X^{m}\rightarrow X$ and $n$-ary
multiplication $\mathbf{\mu}_{n}:X^{n}\rightarrow X$, such that additively
$\left\langle X\mid\mathbf{\nu}_{m}\mid assoc\mid comm\mid solv\right\rangle $
is an $m$-ary commutative group, and multiplicatively $\left\langle
X\mid\mathbf{\mu}_{n}\mid assoc\right\rangle $ is a $n$-ary semigroup, while
the operations $\mathbf{\nu}_{m}$ and $\mathbf{\mu}_{n}$ are not arbitrary,
but connected by the following distributivity property having $n$ relations%
\begin{align}
&  \mathbf{\mu}_{n}\left[  \mathbf{\nu}_{m}\left[  x_{1},\ldots x_{m}\right]
,y_{2},y_{3},\ldots y_{n}\right] \nonumber\\
&  =\mathbf{\nu}_{m}\left[  \mathbf{\mu}_{n}\left[  x_{1},y_{2},y_{3},\ldots
y_{n}\right]  ,\mathbf{\mu}_{n}\left[  x_{2},y_{2},y_{3},\ldots y_{n}\right]
,\ldots\mathbf{\mu}_{n}\left[  x_{m},y_{2},y_{3},\ldots y_{n}\right]  \right]
\label{dis1}\\
&  \mathbf{\mu}_{n}\left[  y_{1},\mathbf{\nu}_{m}\left[  x_{1},\ldots
x_{m}\right]  ,y_{3},\ldots y_{n}\right] \nonumber\\
&  =\mathbf{\nu}_{m}\left[  \mathbf{\mu}_{n}\left[  y_{1},x_{1},y_{3},\ldots
y_{n}\right]  ,\mathbf{\mu}_{n}\left[  y_{1},x_{2},y_{3},\ldots y_{n}\right]
,\ldots\mathbf{\mu}_{n}\left[  y_{1},x_{m},y_{3},\ldots y_{n}\right]  \right]
\label{dis2}\\
&  \vdots\nonumber\\
&  \mathbf{\mu}_{n}\left[  y_{1},y_{2},\ldots y_{n-1},\mathbf{\nu}_{m}\left[
x_{1},\ldots x_{m}\right]  \right] \nonumber\\
&  =\mathbf{\nu}_{m}\left[  \mathbf{\mu}_{n}\left[  y_{1},y_{2},\ldots
y_{n-1},x_{1}\right]  ,\mathbf{\mu}_{n}\left[  y_{1},y_{2},\ldots
y_{n-1},x_{2}\right]  ,\ldots\mathbf{\mu}_{n}\left[  y_{1},y_{2},\ldots
y_{n-1},x_{m}\right]  \right]  , \label{dis3}%
\end{align}
where $x_{i},y_{j}\in X$, $i=1,\ldots,m$, $j=1,\ldots,n$.

If not all relations of distributivity (\ref{dis1})--(\ref{dis3}), as well as
associativity (\ref{mi}) are not satisfied, they are called partial (as
opposed to the total ones), which gives enormous different versions of
polyadic rings.

For instance, the exotic $\left(  3,2\right)  $-ring over rationals
$\mathcal{R}_{3,2}=\left\langle \mathbb{R\mid}\mathbf{\nu}_{3},\mathbf{\mu
}_{2}\right\rangle $ with the derived ternary addition $\mathbf{\nu}%
_{3}\left[  x,y,z\right]  =x\cdot y\cdot z$ and the binary multiplication
$\mathbf{\mu}_{2}\left[  x,y\right]  =x^{y}$ is noncommutative, because
$x^{y}\neq y^{x}$, nonassociative, since $\left(  x^{y}\right)  ^{z}\neq
x^{\left(  y^{z}\right)  }$, and partial (left) distributive, because only the
first relation (\ref{dis1}) holds $\left(  xyz\right)  ^{t}=x^{t}y^{t}z^{t}$,
but $t^{xyz}\neq t^{x}t^{y}t^{z}$, $x,y,z,t\in\mathbb{R}$.

The simplest ternary ring is the commutative, totally associative and totally
distributive $\left(  2,3\right)  $-ring of negative integer numbers
$\mathcal{R}_{2,3}^{\left(  -\right)  }=\left\langle \mathbb{Z}^{-}%
\mathbb{\mid}\mathbf{\nu}_{2},\mathbf{\mu}_{3}\right\rangle $ with the derived
ternary addition $\mathbf{\nu}_{2}\left[  x,y\right]  =x+y$ and the nonderived
ternary multiplication $\mathbf{\mu}_{2}\left[  x,y\right]  =x\cdot y\cdot z$
(the binary product $x\cdot y$ gives result out of the underlying set
$\mathbb{Z}^{-}$ that is in positive integers $\mathbb{Z}^{+}$).

More general polyadic number rings which are formed by representatives of
congruence (residue) classes were introduced in \cite{dup2017a,dup2019}.
Recall that a congruence class of a non-negative integer $a\in\mathbb{Z}^{+}$,
modulo natural $b\in\mathbb{N}$, is defined by%
\begin{equation}
\left[  \left[  a\right]  \right]  _{b}=\left\{  x_{k}^{\left[  a,b\right]
}\mid k\in\mathbb{Z}\right\}  , \label{ab}%
\end{equation}
where%
\begin{align}
x_{k}^{\left[  a,b\right]  }  &  =a+bk,\label{xa}\\
0  &  \leq a\leq b-1, \label{xa1}%
\end{align}
is a representative of the class being infinite. The set of representatives
$X^{\left[  a,b\right]  }\equiv\left\{  x_{k}^{\left[  a,b\right]  }\right\}
$ (as elements of the congruence class $\left[  \left[  a\right]  \right]
_{b}$) was never considered before in details, because it does not form any
binary structure, while the \textquotedblleft internal\textquotedblright%
\ operations are not simultaneously closed $x_{k_{1}}^{\left[  a,b\right]
}+x_{k_{2}}^{\left[  a,b\right]  }\notin X^{\left[  a,b\right]  }$, $x_{k_{1}%
}^{\left[  a,b\right]  }\cdot x_{k_{2}}^{\left[  a,b\right]  }\notin
X^{\left[  a,b\right]  }$ for arbitrary $a$, $b$ (obviously, $X^{\left[
01\right]  }\in\mathbb{Z}$ are ordinary binary integers). This was a reason,
why only operations between classes (the binary class addition $+^{\prime}$
and binary class multiplication $\cdot^{\prime}$) were considered, we call
them \textquotedblleft external\textquotedblright\ operations \cite{dup2017a}.

Nevertheless, for special values of the congruence class parameters $a$ and
$b$ polyadic (nonbinary) operations on $X^{\left[  a,b\right]  }$ can be
defined, and they are closed \cite{dup2017a}. Indeed, the closeness of $m$
additions and $n$ multiplications in $\mathbb{Z}$%
\begin{align}
x_{k_{1}}^{\left[  a,b\right]  }+\ldots+x_{k_{m}}^{\left[  a,b\right]  }  &
\in X^{\left[  a,b\right]  }\subset\mathbb{Z},\label{x1}\\
x_{k_{1}}^{\left[  a,b\right]  }\cdot\ldots\cdot x_{k_{n}}^{\left[
a,b\right]  }  &  \in X^{\left[  a,b\right]  }\subset\mathbb{Z}%
,\ \ \ \ \ \ \ k_{i},m,n\in\mathbb{N}, \label{x2}%
\end{align}
hold valid, if the following \textquotedblleft quantization\textquotedblright%
\ conditions%
\begin{align}
ma^{\prime}  &  \equiv a^{\prime}\left(  \operatorname{mod}b^{\prime}\right)
\Longleftrightarrow\dfrac{\left(  m-1\right)  a^{\prime}}{b^{\prime}%
}=I^{\left(  m\right)  }\left(  a^{\prime},b^{\prime}\right)  \equiv
I=\operatorname{integer},\label{maa}\\
a^{\prime n}  &  \equiv a^{\prime}\left(  \operatorname{mod}b^{\prime}\right)
\Longleftrightarrow\dfrac{a^{\prime n}-a^{\prime}}{b^{\prime}}=J^{\left(
n\right)  }\left(  a^{\prime},b^{\prime}\right)  =J=\operatorname{integer},
\label{ana}%
\end{align}
are satisfied with the special values of parameters as solutions $\left\{
a^{\prime}\right\}  \subset\left\{  a\right\}  $, $\left\{  b^{\prime
}\right\}  \subset\left\{  b\right\}  $, $X^{\left[  a^{\prime},b^{\prime
}\right]  }\subset X^{\left[  a,b\right]  }$. The corresponding mapping is
called the arity shape%
\begin{equation}
\Psi_{m,n}^{\left[  a^{\prime}b^{\prime}\right]  }:\left(  a^{\prime
},b^{\prime}\right)  \longrightarrow\left(  m,n\right)  , \label{as}%
\end{equation}
and it is presented for lowest values in \textsc{Table 1} of \cite{dup2019},
however $\Psi_{m,n}^{\left[  a^{\prime}b^{\prime}\right]  }$ cannot be
presented by a formula. The arity shape mapping (\ref{as}) is injective and
non-surjective (empty cells in \textsc{Table 1}): for the congruence classes
$\left[  \left[  2\right]  \right]  _{4}$, $\left[  \left[  2\right]  \right]
_{8}$, $\left[  \left[  3\right]  \right]  _{9}$, $\left[  \left[  4\right]
\right]  _{8}$, $\left[  \left[  6\right]  \right]  _{8}$ and $\left[  \left[
6\right]  \right]  _{9}$ there are no solutions of the \textquotedblleft
quantization\textquotedblright\ conditions (\ref{maa})--(\ref{ana}), while,
e.g., $m=5$, $n=6$ corresponds to different congruence classes $\left[
\left[  2\right]  \right]  _{5}$, $\left[  \left[  3\right]  \right]  _{5}$,
$\left[  \left[  2\right]  \right]  _{10}$, and $\left[  \left[  8\right]
\right]  _{10}$.

The closeness (\ref{x1})--(\ref{x2}) and the \textquotedblleft
quantization\textquotedblright\ conditions (\ref{maa})--(\ref{ana}) allow us
to define two abstract polyadic operations on the set of the abstract elements
$\mathsf{X}^{\left[  a^{\prime},b^{\prime}\right]  }=\left\{  \mathsf{x}%
_{k}^{\left[  a^{\prime},b^{\prime}\right]  }\right\}  \equiv\left\{
\mathsf{x}_{k}\right\}  $ reflecting the representatives of the fixed
congruence class $\left[  \left[  a^{\prime}\right]  \right]  _{b^{\prime}}$:
$m$-ary addition and $n$-ary multiplication%
\begin{align}
\mathbf{\nu}_{m}\left[  \mathsf{x}_{k_{1}},\ldots,\mathsf{x}_{k_{m}}\right]
&  =\mathsf{x}_{k_{1}}+\ldots+\mathsf{x}_{k_{m}}=\mathsf{x}_{k_{0}}%
,\label{nm}\\
\mathbf{\mu}_{n}\left[  \mathsf{x}_{r_{1}},\ldots,\mathsf{x}_{r_{n}}\right]
&  =\mathsf{x}_{r_{1}}\cdot\ldots\cdot\mathsf{x}_{r_{n}}=\mathsf{x}_{r_{0}},
\label{mn}%
\end{align}
where%
\begin{align}
k_{0}  &  =k_{1}+\ldots+k_{m}+I^{\left(  m\right)  }\left(  a^{\prime
},b^{\prime}\right)  ,\label{k0}\\
r_{0}  &  =s\left(  r_{i},a^{\prime},b^{\prime}\right)  +J^{\left(  n\right)
}\left(  a^{\prime},b^{\prime}\right)  , \label{r0}%
\end{align}
where the invariants $I^{\left(  m\right)  }\left(  a^{\prime},b^{\prime
}\right)  $ and $J^{\left(  n\right)  }\left(  a^{\prime},b^{\prime}\right)  $
are in (\ref{maa})--(\ref{ana}), and the integer $s\left(  r_{i},a^{\prime
},b^{\prime}\right)  $ is a special polynomial of $r_{i},a^{\prime},b^{\prime
}$, which follows from the presentation (\ref{xa}). For instance, in the
ternary case $n=3$ we have $s\left(  r_{1},r_{2},r_{3},a^{\prime},b^{\prime
}\right)  =a^{\prime2}r_{1}+a^{\prime2}r_{2}+a^{\prime2}r_{3}+a^{\prime
}b^{\prime}r_{1}r_{2}+a^{\prime}b^{\prime}r_{1}r_{3}+a^{\prime}b^{\prime}%
r_{2}r_{3}+b^{\prime2}r_{1}r_{2}r_{3}$.

Because the abstract polyadic operations $\mathbf{\nu}_{m}\mathbf{\ }$and
$\mathbf{\mu}_{n}$ are closed, commutative, totally associative and totally
distributive (two latter ones follow from the binary associativity and
distributivity), we can define the abstract commutative polyadic $\left(
m,n\right)  $-ring of integers%
\begin{equation}
\mathrm{Z}_{m,n}^{\left[  a^{\prime},b^{\prime}\right]  }=\left\langle
\mathsf{X}^{\left[  a^{\prime},b^{\prime}\right]  }\mid\mathbf{\nu}%
_{m},\mathbf{\mu}_{n}\mid comm\right\rangle . \label{za}%
\end{equation}

It follows from (\ref{wl}), that we have \textquotedblleft double
quantization\textquotedblright: in the $\left(  m,n\right)  $-ring we can add%
\begin{equation}
w_{\nu}^{admiss}\left(  m\right)  =\ell_{\nu}\left(  n-1\right)  +1
\label{wa1}%
\end{equation}
elements ($\ell_{\nu}$ is quantity of composed $m$-ary additions) and multiply%
\begin{equation}
w_{\mu}^{admiss}\left(  n\right)  =\ell_{\mu}\left(  n-1\right)  +1
\label{wa2}%
\end{equation}
elements, where $\ell_{\mu}$ is quantity of composed $n$-ary multiplications,
to be in the same underlying set $\mathsf{X}^{\left[  a^{\prime},b^{\prime
}\right]  }$. It is obvious, that $\left(  2,2\right)  $-ring $\mathrm{Z}%
_{2,2}^{\left[  0,1\right]  }=\mathbb{Z}$ is the binary ring of ordinary
integers, and $w_{\nu}^{admiss}\left(  m\right)  =\ell_{\nu}$, $w_{\mu
}^{admiss}\left(  n\right)  =\ell_{\mu}$ are any natural numbers $\mathbb{N}$
(without any \textquotedblleft quantization\textquotedblright).

For example, in the polyadic $\left(  8,7\right)  $-ring $\mathrm{Z}%
_{8,7}=\mathrm{Z}_{8,7}^{\left[  5,7\right]  }$ we can add only $7\ell_{\nu
}+1=8,15,22\ldots$ elements $\mathsf{x}_{i}$ and multiply $6\ell_{\mu
}+1=7,13,19\ldots$ ones.

The elements of the $\left(  m,n\right)  $-ring $\mathrm{Z}_{m,n}%
=\mathrm{Z}_{m,n}^{\left[  a^{\prime},b^{\prime}\right]  }$ (\ref{za}) are
abstract variables $\mathsf{x}_{k}$ obeying the $m$-ary addition (\ref{nm})
and $n$-ary multiplication (\ref{mn}) which inherit the \textquotedblleft
internal\textquotedblright\ operations (\ref{x1})--(\ref{x2}) in the
congruence class $\left[  \left[  a^{\prime}\right]  \right]  _{b^{\prime}}$.
Therefore, the elements of $\left(  m,n\right)  $-ring $\mathrm{Z}_{m,n}$
should carry the arities as additional indices to distinguish elements
corresponding to the same representative (as decimal number) of the initial
congruence classes.

For instance, consider two infinite zeroless and unitless nonderived abstract
polyadic rings%
\begin{align}
\mathrm{Z}_{4,3}  &  =\left\{  \ldots\mathsf{-13,-10,-7_{4,3}%
,-4,-1,2,5,8_{4,3},11,14,17}\ldots\right\}  ,\label{z43}\\
\mathrm{Z}_{6,5}  &  =\left\{  \ldots\mathsf{-22,-17,-12,-7_{6,5}%
,-2,3,8_{6,5},13,18,23,28}\ldots\right\}  , \label{z65}%
\end{align}
which are constructed from the congruence classes $\left[  \left[  2\right]
_{3}\right]  $ and $\left[  \left[  3\right]  \right]  _{5}$, respectively.
The intersection of the generating classes for $\left\vert k\right\vert \leq5$
consists of two integer numbers $U_{\left[  3,5\right]  }^{\left[  2,3\right]
}=\left[  \left[  2\right]  _{3}\right]  \cap\left[  \left[  3\right]
\right]  _{5}=\left\{  -7,8\right\}  $ (in general, $U_{\left[  3,5\right]
}^{\left[  2,3\right]  }$ this is an infinite set). We treat the elements of
the rings (\ref{z43})--(\ref{z65}) not as ordinary decimal numbers, but as
abstract variables (or abstract symbols, as, e.g., the letters $A,B,C$ in the
hexadecimal numeral system) $\mathsf{x}_{k}^{\left[  2,3\right]  }$ and
$\mathsf{x}_{k}^{\left[  3,5\right]  }$ obeying (\ref{nm})--(\ref{r0}) and
carrying additional arity lower indices (which are written manifestly in the
needed cases only). In this way, we conclude that the elements corresponding
to the same numbers in the class intersection $U_{\left[  3,5\right]
}^{\left[  2,3\right]  }$ are different $\left(  \mathsf{-7_{4,3}}\right)
\neq\left(  \mathsf{-7_{4,3}}\right)  $, $\mathsf{8_{4,3}}$ $\neq
\mathsf{8_{6,5}}$, in the sense that they obey different operations and their
arities in distinct polyadic rings $\mathrm{Z}_{4,3}$ and $\mathrm{Z}_{6,5}$.
Indeed, e.g., $\mathsf{8_{4,3}}=\mathsf{x}_{2}^{\left[  2,3\right]  }%
\in\mathrm{Z}_{4,3}$ and $\mathsf{8_{6,5}}=\mathsf{x}_{1}^{\left[  3,5\right]
}\in\mathrm{Z}_{6,5}$, and their first polyadic powers (\ref{xl}) are
different $\left(  \mathsf{8_{4,3}}\right)  ^{\left\langle 1\right\rangle
}=\mathbf{\mu}_{3}\left[  \left(  \mathsf{8_{4,3}}\right)  ^{3}\right]
=\mathsf{512}_{4,3}=\mathsf{x}_{102}^{\left[  2,3\right]  }$, while $\left(
\mathsf{8_{6,5}}\right)  ^{\left\langle 1\right\rangle }=\mathbf{\mu}%
_{5}\left[  \left(  \mathsf{8_{6,5}}\right)  ^{5}\right]  =\mathsf{32768}%
_{6,5}=\mathsf{x}_{6553}^{\left[  3,5\right]  }$.

\section{\textsc{Polyadization of positional numeral systems}}

Let us remind the standard positional numeral system, as the presentation of a
number by a special sequence of its digits (for numerous extended versions and
history, see, e.g., \cite{ifrah}).

In the manifest form, the presentation of a number over the (binary) ring of
non-negative integers $\mathbb{Z}^{+}=0,1,2\ldots$ is defined by the base
(radix) $p\in\mathbb{Z}^{+}$ and the digits $y\left(  i\right)  \in
\mathbb{Z}^{+}$ as follows (place-value notation)%
\begin{align}
N^{\left(  \ell\right)  }\left(  p\right)   &  =%
%TCIMACRO{\tsum \limits_{i=0}^{\ell-1}}%
%BeginExpansion
{\textstyle\sum\limits_{i=0}^{\ell-1}}
%EndExpansion
y\left(  i\right)  p^{i}=y\left(  \ell-1\right)  p^{\ell-1}+y\left(
\ell-2\right)  p^{\ell-2}+\ldots+y\left(  1\right)  p+y\left(  0\right)
\nonumber\\
&  \Longrightarrow\left(  y\left(  \ell-1\right)  y\left(  \ell-2\right)
\ldots y\left(  1\right)  y\left(  0\right)  \right)  _{p}^{\left(
\ell\right)  },\label{n}\\
0 &  \leq y\left(  i\right)  \leq p-1,\label{n1}%
\end{align}
where the natural $\ell=l\left(  p\right)  \in\mathbb{N}$ is the quantity of
digits (function of the base $p$), being simultaneously the amount of terms in
the l.h.s., $\ell-1$ (in the binary case only) coincides with the quantity of
additions and the total amount of multiplications in the first term, and
usually $p^{0}=1$ (to have the same summation formula with $i=0$). Commonly,
the leading zeroes in the r.h.s. of (\ref{n}) are omitted. The change of the
base $p\mapsto p^{\prime}$ leads to another presentations of the same number,
usually with the different quantity of digits $\ell\mapsto\ell^{\prime
}=l\left(  p^{\prime}\right)  $%
\begin{equation}
N^{\left(  \ell^{\prime}\right)  }\left(  p^{\prime}\right)  =%
%TCIMACRO{\tsum \limits_{i=0}^{\ell^{\prime}-1}}%
%BeginExpansion
{\textstyle\sum\limits_{i=0}^{\ell^{\prime}-1}}
%EndExpansion
y^{\prime}\left(  i\right)  p^{\prime i}\Longrightarrow\left(  y^{\prime
}\left(  \ell-1\right)  y^{\prime}\left(  \ell-2\right)  \ldots y^{\prime
}\left(  1\right)  y^{\prime}\left(  0\right)  \right)  _{p^{\prime}}^{\left(
\ell^{\prime}\right)  },\ \ 0\leq y^{\prime}\left(  i\right)  \leq p^{\prime
}-1.
\end{equation}

In the binary case, using $\ell$ digits one can describe $E\left(  p\right)
=p^{\ell}$ numbers: $0,\ldots p^{\ell}-1\in\mathbb{Z}^{+}$. The efficiency of
a numeral system is the ability to represent as many numbers as possible using
the smallest total number of symbols $s$. In this case, the number of digits
becomes $s/p$, and the quantity of the described numbers (the efficiency
function) is%
\begin{equation}
E\left(  p\right)  =p^{\frac{s}{p}}. \label{e}%
\end{equation}

The function $E\left(  p\right)  $ reaches its maximum (in $p$), when
$p=e\approx2.718...$ (the Euler's number), and therefore the most efficient
numeral system (with integer $p$) is $p_{\max}=3$.

For further details about positional numeral systems (over the binary ring
$\mathbb{Z}^{+}$), see, e.g., \cite{ifrah}, and refs. therein.

Now we generalize the above construction to the polyadic integer numbers
\cite{dup2017a,dup2019}, that is we build an analog of positional numeral
system over polyadic $\left(  m,n\right)  $-rings $\mathrm{Z}_{m,n}$
considered in the previous section. The polyadic analog of binary numeral
systems have complicated and nontrivial structure, because the underlying
$\left(  m,n\right)  $-rings obey unusual peculiarities, for instance, some of
them can not have zero or/and unity, absence of the natural ordering, etc.
(see, for details, \cite{duplij2022}).

First, we express the binary positional numeral system (presentation
(\ref{n})) in the polyadic notation from \textsc{Section} \ref{sec-prel}.
Indeed, let the binary ring of integers is $\mathbb{Z}=\left\langle \nu
,\mu\mid comm\mid assoc\right\rangle $, where $\nu$ and $\mu$ are the ordinary
binary addition $\nu\left[  x,y\right]  =x+y$ and multiplication $\mu\left[
x,y\right]  =y\cdot y$, $x,y\in\mathbb{Z}$. In this notation, each term in
(\ref{n}) can be written as the composition of $i+1$ multiplications (see
(\ref{ml})) $\mu^{\circ\left(  i+1\right)  }\left[  y\left(  i\right)
,\overset{i}{\overbrace{p,\ldots,p}}\right]  \equiv\mu^{\circ\left(
i+1\right)  }\left[  y\left(  i\right)  ,p^{i}\right]  $, the first term
becomes $\mu^{\circ\left(  \ell_{\mu}-1\right)  }\left[  y\left(  \ell_{\mu
}-1\right)  ,p^{\ell_{\mu}-1}\right]  $, and the sum as the composition of
$\ell_{\nu}$ additions, such that the place-value presentation takes the form
of operation compositions%
\begin{align}
N^{\left(  \ell_{\nu},\ell_{\mu}\right)  }\left(  p\right)   &  =\nu
^{\circ\ell_{\nu}}\left[  \overset{\ell_{\nu}+1}{\overbrace{\mu^{\circ\left(
\ell_{\mu}-1\right)  }\left[  y\left(  \ell_{\mu}-1\right)  ,p^{\ell_{\mu}%
-1}\right]  ,\mu^{\circ\left(  \ell_{\mu}-1\right)  }\left[  y\left(
\ell_{\mu}-2\right)  ,p^{\ell_{\mu}-2}\right]  ,\ldots,\mu\left[  y\left(
1\right)  ,p\right]  ,y\left(  0\right)  }}\right]  \nonumber\\
&  \Longrightarrow\left(  \overset{\ell_{\nu}+1}{\overbrace{y\left(  \ell
_{\mu}-1\right)  y\left(  \ell_{\mu}-2\right)  \ldots y\left(  1\right)
y\left(  0\right)  }}\right)  _{p}^{\left(  \ell_{v},\ell_{\mu}\right)
},\ \ \ 0\leq y\left(  i\right)  \leq p-1\in\mathbb{Z}^{+},\label{nl}%
\end{align}
where the important additional consistency condition%
\begin{equation}
\ell_{\mu}=\ell_{\nu}+1=\ell,\label{ll}%
\end{equation}
which follows from the construction itself. Note that in the composition form
we do not write multiplier for the last degit $y\left(  0\right)  $ at all for
consistency with the higher arity generalizations (not all multiplicative
parts of polyadic rings contain unity, see \textsc{Section \ref{sec-polring}},
however we need the last digit in any case).

For example, in the case $p=7$ and $\ell=\ell_{\mu}=\ell_{\nu}+1=3$ we have
the ordinary binary positional numeral presentation in the polyadic notation%
\begin{equation}
\left(  165\right)  _{7}^{\left(  2,2\right)  }=\nu^{\circ2}\left[  \mu
^{\circ2}\left[  1,7,7\right]  ,\mu\left[  6,7\right]  ,5\right]
=1\cdot7\cdot7+6\cdot7+5=\left(  96\right)  _{10}^{\left(  1,1\right)  }.
\end{equation}

Now the polyadization of the binary place-value presentation in the
composition form (\ref{nl}) can be done in the straightforward way. We assume
that all ordinary numbers become polyadic numbers, that is the initial number
ring should be exchange to the polyadic $\left(  m,n\right)  $-ring:
$\mathbb{Z=}\left\langle \nu,\mu\right\rangle \longrightarrow\mathrm{Z}%
_{m,n}=\left\langle \mathbf{\nu}_{m},\mathbf{\mu}_{n}\right\rangle $. Thus,
modifying (\ref{nl}) consistently by $\nu\rightarrow\mathbf{\nu}_{m}$ and
$\mu\rightarrow\mathbf{\mu}_{n}$, taking into account the admissible length of
words (\ref{wa1})-(\ref{wa2}), we propose the direct polyadic generalization
of the standard place-value presentation (\ref{n}) of a polyadic number
$\mathbf{N}\in\mathrm{Z}_{m,n}$ by the polyadic digits $\mathbf{y}\left(
i\right)  \in\mathrm{Z}_{m,n}$ and the polyadic base $\mathbf{p}\in
\mathrm{Z}_{m,n}$%
\begin{align}
&  \mathbf{N}_{m,n}^{\left(  \ell_{\nu},\ell_{\mu}\right)  }\left(
\mathbf{p}\right)  =\mathbf{\nu}_{m}^{\circ\ell_{\nu}}\left[  {}\right.
\nonumber\\
&  \left.  \overset{\ell_{\nu}\left(  m-1\right)  +1}{\overbrace{\mathbf{\mu
}_{n}^{\circ\left(  \ell_{\mu}-1\right)  }\left[  \mathbf{y}\left(  \ell_{\mu
}-1\right)  ,\mathbf{p}^{\ell_{\mu}\left(  n-1\right)  }\right]  ,\mathbf{\mu
}_{n}^{\circ\left(  \ell_{\mu}-1\right)  }\left[  \mathbf{y}\left(  \ell_{\mu
}-2\right)  ,\mathbf{p}^{\ell_{\mu}\left(  n-1\right)  -1}\right]
,\ldots,\mathbf{\mu}_{n}\left[  \mathbf{y}\left(  1\right)  ,\mathbf{p}%
^{n-1}\right]  ,\mathbf{y}\left(  0\right)  }}\right]  \nonumber\\
&  \Longrightarrow\left(  \overset{\ell_{\nu}\left(  m-1\right)
+1}{\overbrace{\mathbf{y}\left(  \ell_{\mu}-1\right)  \mathbf{y}\left(
\ell_{\mu}-2\right)  \ldots\mathbf{y}\left(  1\right)  \mathbf{y}\left(
0\right)  }}\right)  _{m,n;\mathbf{p}}^{\left(  \ell_{v},\ell_{\mu}\right)
},\label{nn}%
\end{align}
where the arity indices $m,n$, when they are obvious, can be omitted for
conciseness. The polyadic analog of the consistency condition (\ref{ll}) now
is%
\begin{equation}
\ell_{\mu}=\ell_{\nu}\left(  m-1\right)  +1,\label{lm}%
\end{equation}
and therefore the minimal quantity of multiplications in (\ref{nn}) is
$\ell_{\mu}\geq m$, which coincides with minimal amount of digits in the
place-value presentation over $\left(  m,n\right)  $-ring.

\begin{theorem}
In the polyadic numeral system over $\left(  m,n\right)  $-ring the minimal
number of digits is more or equal than the arity of addition $m$.
\end{theorem}

The $m$-ary addition of polyadic numbers $\mathbf{N}$ in the place-value
presentation (\ref{nn}) can be done using the total polyadic distributivity
(\ref{dis1})--(\ref{dis3}) in the general form by adding digits (as in the
binary case, prime sign is not a derivative, but numerates variables)%
\begin{align}
\mathbf{\nu}_{m}\left[  \mathbf{N}_{m,n}^{\prime\left(  \ell_{\nu},\ell_{\mu
}\right)  }\left(  \mathbf{p}\right)  ,\mathbf{N}_{m,n}^{\prime\prime\left(
\ell_{\nu},\ell_{\mu}\right)  }\left(  \mathbf{p}\right)  ,\ldots
,\mathbf{N}_{m,n}^{\prime\prime\prime\left(  \ell_{\nu},\ell_{\mu}\right)
}\left(  \mathbf{p}\right)  \right]   &  =\mathbf{N}_{m,n}^{\left(  \ell_{\nu
},\ell_{\mu}\right)  }\left(  \mathbf{p}\right)  ,\\
\mathbf{\nu}_{m}\left[  \mathbf{y}^{\prime}\left(  i\right)  ,\mathbf{y}%
^{\prime\prime}\left(  i\right)  ,\ldots,\mathbf{y}^{\prime\prime\prime
}\left(  i\right)  \right]   &  =\mathbf{y}\left(  i\right)  ,\ \ \ i=0,\ldots
,\ell_{\mu}.
\end{align}

The $n$-ary multiplication of polyadic numbers $\mathbf{N}$ is more
complicated and should be made in each concrete case manifestly.

In general, the direct polyadization formula (\ref{nn}) for place-value
presentation can be considered for any commutative polyadic $\left(
m,n\right)  $-ring.

Here we will study some examples for polyadic rings of integer numbers from
the previous \textsc{Section} \textsc{\ref{sec-polring}}. Let us consider the
abstract polyadic number ring $\mathrm{Z}_{4,3}^{\left[  2,3\right]  }$
(\ref{z43}) generated by the congruence class $\left[  \left[  2\right]
\right]  _{3}$. The polyadic integer numbers $\mathsf{x}_{k}^{\left[
2,3\right]  }$ are \textquotedblleft symmetric\textquotedblright\ with respect
to $\mathsf{x}_{k=0}^{\left[  2,3\right]  }=a=2$ (playing the role of zero in
the binary case $a=0$, $b=1$ and $x_{k}=k$), therefore instead of (\ref{xa1}),
we use%
\begin{equation}
0\leq k_{y\left(  i\right)  }\leq k_{p}-1,
\end{equation}
where $\mathbf{y}\left(  i\right)  =\mathsf{x}_{k_{y\left(  i\right)  }}$,
$\mathbf{p}=\mathsf{x}_{k_{p}}$, and $\mathsf{x}_{k}=\mathsf{x}_{k}^{\left[
2,3\right]  }\in\mathsf{X}^{\left[  2,3\right]  }$. In this way, for the
simplest polyadic base $\mathbf{p}=\mathsf{x}_{2}=\mathsf{8}$, $\ell_{\nu}=1$
and $\ell_{\mu}=4$ (see (\ref{lm})) we have the $4$-digit \textbf{(}%
$\mathbf{y}\left(  i\right)  =\mathsf{2,5}$, $i=0,1,2,3$) polyadic numeral
presentation for a number from $\mathrm{Z}_{4,3}^{\left[  2,3\right]
}\mathrm{\ }$as%
\begin{align}
&  \mathbf{N}_{4,3}^{\left(  1,3\right)  }\left(  \mathsf{8}\right)
=\mathbf{\nu}_{4}\left[  \mathbf{\mu}_{3}^{\circ3}\left[  \mathbf{y}\left(
3\right)  ,\mathsf{8},\mathsf{8},\mathsf{8},\mathsf{8},\mathsf{8}%
,\mathsf{8}\right]  ,\mathbf{\mu}_{3}^{\circ2}\left[  \mathbf{y}\left(
2\right)  ,\mathsf{8},\mathsf{8},\mathsf{8},\mathsf{8}\right]  ,\mathbf{\mu
}_{3}\left[  \mathbf{y}\left(  1\right)  ,\mathsf{8},\mathsf{8}\right]
,\mathbf{y}\left(  0\right)  \right]  \nonumber\\
&  \Longrightarrow\left(  \mathbf{y}\left(  3\right)  \mathbf{y}\left(
2\right)  \mathbf{y}\left(  1\right)  \mathbf{y}\left(  0\right)  \right)
_{4,3;\mathsf{8}}^{\left(  1,3\right)  }\equiv\left(  \mathbf{y}\left(
3\right)  \mathbf{y}\left(  2\right)  \mathbf{y}\left(  1\right)
\mathbf{y}\left(  0\right)  \right)  _{\mathsf{8}}.
\end{align}

The polyadic numerals (polyadic numbers of the $\left(  m,n\right)  $-ring
$\mathrm{Z}_{4,3}^{\left[  2,3\right]  }$ that can be presented in
base-$\mathsf{8}$ place-value form by $4$ digits) are%
\begin{align}
\left(  \mathsf{2,2,2,2}\right)  _{\mathsf{8}}  &  =\mathsf{532610}%
=\mathsf{x}_{177536},\left(  \mathsf{2,2,2,5}\right)  _{\mathsf{8}%
}=\mathsf{532613}=\mathsf{x}_{177537},\left(  \mathsf{2,2,5,2}\right)
_{\mathsf{8}}=\mathsf{532802}=\mathsf{x}_{177600},\nonumber\\
\left(  \mathsf{2,2,5,5}\right)  _{\mathsf{8}}  &  =\mathsf{532805}%
=\mathsf{x}_{177601},\left(  \mathsf{2,5,2,2}\right)  _{\mathsf{8}%
}=\mathsf{544898}=\mathsf{x}_{181632},\left(  \mathsf{2,5,2,5}\right)
_{\mathsf{8}}=\mathsf{544901}=\mathsf{x}_{181633},\nonumber\\
\left(  \mathsf{2,5,5,2}\right)  _{\mathsf{8}}  &  =\mathsf{545090}%
=\mathsf{x}_{181696},\left(  \mathsf{2,5,5,5}\right)  _{\mathsf{8}%
}=\mathsf{545093}=\mathsf{x}_{181697},\left(  \mathsf{5,2,2,2}\right)
_{\mathsf{8}}=\mathsf{1319042}=\mathsf{x}_{439680},\nonumber\\
\left(  \mathsf{5,2,2,5}\right)  _{\mathsf{8}}  &  =\mathsf{1319045}%
=\mathsf{x}_{439681},\left(  \mathsf{5,2,5,2}\right)  _{\mathsf{8}%
}=\mathsf{1319234}=\mathsf{x}_{439744},\left(  \mathsf{5,2,5,5}\right)
_{\mathsf{8}}=\mathsf{1319237}=\mathsf{x}_{439745},\nonumber\\
\left(  \mathsf{5,5,2,2}\right)  _{\mathsf{8}}  &  =\mathsf{1331330}%
=\mathsf{x}_{443776},\left(  \mathsf{5,5,2,5}\right)  _{\mathsf{8}%
}=\mathsf{1331333}=\mathsf{x}_{443777},\left(  \mathsf{5,5,5,2}\right)
_{\mathsf{8}}=\mathsf{1331522}=\mathsf{x}_{443840},\nonumber\\
\left(  \mathsf{5,5,5,5}\right)  _{\mathsf{8}}  &  =\mathsf{1331525}%
=\mathsf{x}_{443841}, \label{2}%
\end{align}
which correspond to the ordinary base-$2$ (or binary) numerals with two digits
$y\left(  i\right)  =0,1$.

The more complicated case is the $7$-digit \textbf{(}$\mathbf{y}\left(
i\right)  =\mathsf{2,5,8,11}$, $i=0,\ldots,7$) polyadic numeral presentation
for numbers from $\mathrm{Z}_{4,3}^{\left[  2,3\right]  }\mathrm{\ }$with the
polyadic base (see (\ref{z43})) $\mathbf{p}=\mathsf{x}_{4}=\mathsf{14}$,
composition of two $m$-ary additions $\ell_{\nu}=2$ and seven $n$-ary
multiplications $\ell_{\mu}=7$ (see (\ref{lm}))%
\begin{align}
&  \mathbf{N}_{4,3}^{\left(  1,3\right)  }\left(  \mathsf{14}\right)
=\mathbf{\nu}_{4}^{\circ2}\left[  \mathbf{\mu}_{3}^{\circ6}\left[
\mathbf{y}\left(  6\right)  ,\mathsf{14}^{12}\right]  ,\mathbf{\mu}_{3}%
^{\circ5}\left[  \mathbf{y}\left(  5\right)  ,\mathsf{14}^{10}\right]
,\mathbf{\mu}_{3}^{\circ4}\left[  \mathbf{y}\left(  4\right)  ,\mathsf{14}%
^{8}\right]  ,\mathbf{\mu}_{3}^{\circ3}\left[  \mathbf{y}\left(  3\right)
,\mathsf{14}^{6}\right]  ,\right.  \nonumber\\
&  \left.  \mathbf{\mu}_{3}^{\circ2}\left[  \mathbf{y}\left(  2\right)
,\mathsf{14},\mathsf{14},\mathsf{14},\mathsf{14}\right]  ,\mathbf{\mu}%
_{3}\left[  \mathbf{y}\left(  1\right)  ,\mathsf{14},\mathsf{14}\right]
,\mathbf{y}\left(  0\right)  \right]  \nonumber\\
&  \Longrightarrow\left(  \mathbf{y}\left(  6\right)  \mathbf{y}\left(
5\right)  \mathbf{y}\left(  4\right)  \mathbf{y}\left(  3\right)
\mathbf{y}\left(  2\right)  \mathbf{y}\left(  1\right)  \mathbf{y}\left(
0\right)  \right)  _{4,3;\mathsf{14}}^{\left(  1,3\right)  }\equiv\left(
\mathbf{y}\left(  6\right)  \mathbf{y}\left(  5\right)  \mathbf{y}\left(
4\right)  \mathbf{y}\left(  3\right)  \mathbf{y}\left(  2\right)
\mathbf{y}\left(  1\right)  \mathbf{y}\left(  0\right)  \right)
_{\mathsf{14}}.\label{n5}%
\end{align}

The first several numerals in (\ref{n5}) are%
\begin{align}
\left(  \mathsf{2,2,2,2,2,2,2}\right)  _{\mathsf{14}}  &
=\mathsf{113\ 969\ 300\ 774\ 954},\nonumber\\
\left(  \mathsf{2,2,2,2,2,2,11}\right)  _{\mathsf{14}}  &
=\mathsf{113\ 969\ 300\ 774\ 963},\nonumber\\
\left(  \mathsf{2,2,2,2,2,5,2}\right)  _{\mathsf{14}}  &
=\mathsf{113\ 969\ 300\ 775\ 542},\nonumber\\
\left(  \mathsf{2,2,2,2,2,8,2}\right)  _{\mathsf{14}}  &
=\mathsf{113\ 969\ 300\ 776\ 130},\nonumber\\
\left(  \mathsf{2,2,2,2,2,11,11}\right)  _{\mathsf{14}}  &
=\mathsf{113\ 969\ 300\ 776\ 727,}\nonumber\\
\left(  \mathsf{2,5,2,2,2,2,2}\right)  _{\mathsf{14}}  &
=\mathsf{114\ 837\ 064\ 739\ 882},\nonumber\\
\left(  \mathsf{2,11,2,2,2,2,2}\right)  _{\mathsf{14}}  &
=\mathsf{116\ 572\ 592\ 669\ 738},\nonumber\\
\left(  \mathsf{11,2,2,2,2,2,2}\right)  _{\mathsf{14}}  &
=\mathsf{624\ 214\ 512\ 152\ 618},\nonumber\\
\left(  \mathsf{11,11,11,11,11,11,11}\right)  _{\mathsf{14}}  &
=\mathsf{626\ 831\ 154\ 262\ 247}. \label{14}%
\end{align}

It follows from the examples (\ref{2}) and (\ref{14}), that not all polyadic
integer numbers $\mathbf{N}$ can be represented by the polyadic numeral
formula in the place-value form (\ref{nn}). We call such numbers numerally
representable $\left\{  \mathbf{N}_{reps}\right\}  \subseteq\left\{
\mathbf{N}\right\}  $.

\begin{theorem}
The set of polyadic numbers which are representable in the numeral form is a
subset of the set of all numbers, while the equality is reached in the binary
case only $\mathbf{N}_{reps}=\mathbf{N}\Longleftrightarrow\mathbf{N}%
_{reps},\mathbf{N\in}\mathrm{Z}_{2,2}^{\left[  0,1\right]  }\equiv\mathbb{Z}$.
In the polyadic nonbinary case $m\geq3\vee n\geq3$, the set of representable
numbers is a proper (strict) subset of all numbers $\left\{  \mathbf{N}%
_{reps}\right\}  \subset\left\{  \mathbf{N}\right\}  $.
\end{theorem}

Now we introduce the polyadic analog of the mixed-base (radix) positional
numeral system. In the polyadic notation, the mixed-base version of the binary
(\ref{nl}) can be written as the sum as the composition of $\ell_{\nu}$
additions having different bases $p\left(  i\right)  $ for each digit
$y\left(  i\right)  $, and we assume that for the last digit $y\left(
0\right)  $ the basis multiplier is absent (see note after (\ref{ll})). In
such assumption the connection between the quantity of digits $N_{y}$ and the
quantity of bases $N_{p}$ is%
\begin{equation}
N_{p}=N_{y}-1.\label{np}%
\end{equation}

Then the place-value presentation takes the form of operation compositions%
\begin{align}
N^{\left(  \ell_{\nu},\ell_{\mu}\right)  }\left(  p\right)   &  =\nu
^{\circ\ell_{\nu}}\left[  {}\right.  \nonumber\\
&  \left.  \overset{\ell_{\nu}+1}{\overbrace{\mu^{\circ\left(  \ell_{\mu
}-1\right)  }\left[  y\left(  \ell_{\mu}-1\right)  ,p\left(  \ell_{\mu
}-1\right)  \right]  ,\mu^{\circ\left(  \ell_{\mu}-1\right)  }\left[  y\left(
\ell_{\mu}-2\right)  ,p\left(  \ell_{\mu}-2\right)  \right]  ,\ldots
,\mu\left[  y\left(  1\right)  ,p\left(  1\right)  \right]  ,y\left(
0\right)  }}\right]  \nonumber\\
&  \Longrightarrow\left(  \overset{\ell_{\nu}+1}{\overbrace{y\left(  \ell
_{\mu}-1\right)  y\left(  \ell_{\mu}-2\right)  \ldots y\left(  1\right)
y\left(  0\right)  }}\right)  _{p\left(  \ell_{\mu}-1\right)  ,p\left(
\ell_{\mu}-2\right)  \ldots p\left(  1\right)  }^{\left(  \ell_{v},\ell_{\mu
}\right)  },\ 0\leq y\left(  i\right)  \leq p\left(  i\right)  -1\in
\mathbb{Z}^{+}.\label{mix}%
\end{align}

The simplest binary example of the mixed-base numeral system is computation of
currency in different banknote numbers $y\left(  i\right)  $ of the
denomination $p\left(  i\right)  $. Sometimes the bases are composed and
connected between them taking the \textquotedblleft
recurrent\textquotedblright\ form%
\begin{equation}
p\left(  i+1\right)  =p\left(  i\right)  p\left(  i-1\right)  \ldots p\left(
1\right)  .\label{pi}%
\end{equation}
A binary example of mixed-base recurrent numeral system is the time
calculation in seconds (\textquotedblleft clock\textquotedblright%
/\textquotedblleft timer\textquotedblright) using number of days $y\left(
days\right)  $, hours $y\left(  hours\right)  $, minutes $y\left(
\text{\textit{min}}\right)  $, seconds $y\left(  \text{\textit{seconds}%
}\right)  $, and the lengths $p\left(  days\right)  =24$, $p\left(
hours\right)  =60$, $p\left(  \text{\textit{minutes}}\right)  =60$. E.g. the
polyadic/composition form (\ref{mix}) of \textquotedblleft2 days, 7 hours, 35
minutes, 48 seconds\textquotedblright\ (having $m=n=2$, $\ell_{\nu}=2$,
$\ell_{\mu}=4$), is (in seconds)%
\begin{equation}
\nu^{\circ2}\left[  \mu^{\circ3}\left[  2,24,60,60\right]  ,\mu^{\circ
2}\left[  7,60,60\right]  ,\mu\left[  35,60\right]  ,48\right]  =2\cdot
24\cdot60\cdot60+7\cdot60\cdot60+35\cdot60+48=200\,148.\label{cl}%
\end{equation}

In general, the polyadic analog of the mixed-base (radix) positional numeral
system can be defined by the formula similar to (\ref{nn}), but the base in
each term in the $m$-ary sum will depend of its place $\mathbf{p}%
\longrightarrow\mathbf{p}\left(  i\right)  $, that is%
\begin{align}
&  \mathbf{N}_{m,n}^{\left(  \ell_{\nu},\ell_{\mu}\right)  }\left(
\mathbf{p}\left(  1\right)  ,\ldots\mathbf{p}\left(  \ell_{\mu}-1\right)
\right)  =\mathbf{\nu}_{m}^{\circ\ell_{\nu}}\left[  {}\right.  \nonumber\\
&  \mathbf{\mu}_{n}^{\circ\left(  \ell_{\mu}-1\right)  }\left[  \mathbf{y}%
\left(  \ell_{\mu}-1\right)  ,\mathbf{p}\left(  \ell_{\mu}-1\right)
^{\ell_{\mu}\left(  n-1\right)  }\right]  ,\mathbf{\mu}_{n}^{\circ\left(
\ell_{\mu}-2\right)  }\left[  \mathbf{y}\left(  \ell_{\mu}-2\right)
,\mathbf{p}\left(  \ell_{\mu}-2\right)  ^{\ell_{\mu}\left(  n-1\right)
-1}\right]  ,\ldots\nonumber\\
&  \left.  \ldots,\mathbf{\mu}_{n}\left[  \mathbf{y}\left(  1\right)
,\mathbf{p}\left(  1\right)  ^{n-1}\right]  ,\mathbf{y}\left(  0\right)
\right]  \Longrightarrow\left(  \overset{\ell_{\nu}\left(  m-1\right)
+1}{\overbrace{\mathbf{y}\left(  \ell_{\mu}-1\right)  \mathbf{y}\left(
\ell_{\mu}-2\right)  \ldots\mathbf{y}\left(  1\right)  \mathbf{y}\left(
0\right)  }}\right)  _{m,n;\mathbf{p}\left(  1\right)  ,\ldots\mathbf{p}%
\left(  \ell_{\mu}-1\right)  }^{\left(  \ell_{v},\ell_{\mu}\right)
},\label{mix1}\\
&  \mathbf{p}\left(  1\right)  ,\ldots\mathbf{p}\left(  \ell_{\mu}-1\right)
\in\mathrm{Z}_{m,n}.\nonumber
\end{align}

Here the connection between the quantity of digits $N_{y}$ and the quantity of
bases $N_{p}$ is the same as in the binary case (\ref{np}). Nevertheless, we
can extend each $i$th term using in it different bases (which now have two
indices $\mathbf{p}_{j}\left(  i\right)  \in\mathrm{Z}_{m,n}$) as%
\begin{equation}
\mathbf{\mu}_{n}^{\circ\left(  i-1\right)  }\left[  \mathbf{y}\left(
i-1\right)  ,\mathbf{p}_{1}\left(  i-1\right)  ,\mathbf{p}_{2}\left(
i-1\right)  ,\ldots,\mathbf{p}_{i\left(  n-1\right)  }\left(  i-1\right)
\right]  .\label{miy}%
\end{equation}

Summing up in (\ref{mix1}), instead of (\ref{np}), for the quantity of
polyadic bases $\mathbf{p}_{j}\left(  i\right)  $ through the number of digits
$N_{y}$ we will have%
\begin{equation}
N_{p}=\frac{1}{2}N_{y}\left(  N_{y}-1\right)  \left(  n-1\right)  .
\end{equation}

For example, in the case $\ell_{\nu}=1$ and $\ell_{\mu}=4$ (see (\ref{lm})) we
have the $4$-digit \textbf{(}$\mathbf{y}\left(  i\right)  $, $i=0,1,2,3$)
polyadic mixed-base numeral presentation for a number from the $\left(
4,3\right)  $-ring $\mathrm{Z}_{4,3}^{\left[  2,3\right]  }\mathrm{\ }$as%
\begin{align}
&  \mathbf{N}_{4,3}^{\left(  1,3\right)  }\left(  \mathbf{p}_{j}\left(
i\right)  \right)  =\mathbf{\nu}_{4}\left[  \mathbf{\mu}_{3}^{\circ3}\left[
\mathbf{y}\left(  3\right)  ,\mathbf{p}_{1}\left(  3\right)  ,\mathbf{p}%
_{2}\left(  3\right)  ,\mathbf{p}_{3}\left(  3\right)  ,\mathbf{p}_{4}\left(
3\right)  ,\mathbf{p}_{5}\left(  3\right)  ,\mathbf{p}_{6}\left(  3\right)
\right]  ,\right.  \nonumber\\
&  \left.  \mathbf{\mu}_{3}^{\circ2}\left[  \mathbf{y}\left(  2\right)
,\mathbf{p}_{1}\left(  2\right)  \mathbf{p}_{2}\left(  2\right)
\mathbf{p}_{3}\left(  2\right)  \mathbf{p}_{4}\left(  2\right)  \right]
,\mathbf{\mu}_{3}\left[  \mathbf{y}\left(  1\right)  ,\mathbf{p}_{1}\left(
1\right)  ,\mathbf{p}_{2}\left(  1\right)  \right]  ,\mathbf{y}\left(
0\right)  \right]  \nonumber\\
&  \Longrightarrow\left(  \mathbf{y}\left(  3\right)  \mathbf{y}\left(
2\right)  \mathbf{y}\left(  1\right)  \mathbf{y}\left(  0\right)  \right)
_{4,3;\mathbf{p}_{j}\left(  i\right)  }^{\left(  1,3\right)  }\equiv\left(
\mathbf{y}\left(  3\right)  \mathbf{y}\left(  2\right)  \mathbf{y}\left(
1\right)  \mathbf{y}\left(  0\right)  \right)  _{\mathbf{p}_{j}\left(
i\right)  }.\label{n4p}%
\end{align}
Thus, the mixed-base polyadic $4$-digit numerals over $\left(  4,3\right)
$-ring with $\ell_{\nu}=1$ and $\ell_{\mu}=4$ are described by $12$ different
polyadic bases $\mathbf{p}_{j}\left(  i\right)  \in\mathrm{Z}_{4,3}^{\left[
2,3\right]  }$. Comparing with the \textquotedblleft binary
clock\textquotedblright\ (\ref{mix})--(\ref{cl}), we can call the presented
construction (\ref{mix1})--(\ref{n4p}) as the \textquotedblleft polyadic
clock\textquotedblright.

\section{\textsc{Conclusions}}

In this paper we showed that by transplanting the entire positional-numeral
paradigm into the realm of $\left(  m,n\right)  $-rings we uncovered a
landscape where digit strings, carries, and even the very concept of
\textquotedblleft base\textquotedblright\ behave in ways impossible in the
ordinary arithmetic. Let us distil the main messages.

\textbf{1}. Polyadic rings admit a native positional calculus presented here.
Replacing plus by an $m$-ary sum and product by an $n$-ary multiplication
forces a \textquotedblleft double quantization\textquotedblright\ of word length.

\textbf{2}. Minimal-digit and representability theorems. We proved that in the
positional expansion every numeral must respect the constraint $\ell
_{mult}=\ell_{add}\,\left(  m-1\right)  +1$, and so the shortest admissible
numeral contains more positions (digits) than the additive arity $\ell
_{mult}\geq m$ (Theorem 4.1), and that for $m,n\geq3$ not every element of a
polyadic ring possesses a finite positional expansion (Theorem 4.2). The gap
is governed by the arity-shape invariants $I^{(m)}$ and $J^{(n)}$ originating
in the polyadic features of congruence-class geometry.

\textbf{3}. Mixed bases are generalized naturally yielding a \textquotedblleft
polyadic clock\textquotedblright. Allowing each digit to sit atop its own
$n$-ary tower of bases replaces the linear count $N_{p}=N_{y}-1$ by the
quadratic formula
\[
N_{p}=\frac{1}{2}N_{y}(N_{y}-1)(n-1),
\]
exploding the design space of numeral systems. Our exotic polyadic clock over
$(4,3)$-ring, specified by twelve independent bases, hints at rich
cryptographic and coding-theoretic applications.

\textbf{4}. Concrete catalogues. The worked-out examples in the $(4,3)$ and
$(6,5)$ integer rings illustrate how ordinary integers lift to abstract
polyadic variables that share a decimal "shadow" yet obey distinct algebraic
laws. These catalogues provide ready test-beds for algorithmic experimentation.

Future directions. The present framework raises more questions than it settles:

Algorithms and complexity. What is the cost of addition, \textquotedblleft
carry\textquotedblright\ propagation, and comparison in a polyadic positional
machine? Can we find arities for which arithmetic accelerates relative to
binary RAM models?

Floating-point and analysis. A polyadic analogue of IEEE-754 would require
rounding rules compatible with $m$-ary adders and $n$-ary multipliers. How do
error bounds scale with the arity pair $(m,n)$?

Hardware realization. Ternary logic became tangible with modern CMOS; could
polyadic ALUs (arithmetic logic units) unlock energy or area advantages?

Category-theoretic reformulation. Viewing $(m,n)$-rings as $\mathbb{N}^{2}%
$-graded monoidal categories may clarify the functorial behavior of numeral expansions.

Quantum and cryptographic angles. The inherent nonassociativity of many
polyadic rings could serve as a source of hard problems, while the enlarged
digit alphabet suggests fresh qudit encodings.

\bigskip

%\section{Acknowledgments}

\textbf{Acknowledgments}. The author would like to express his deep
thankfulness to Vladimir Tkach for productive discussions and to Qiang Guo and
Raimund Vogl for valuable support.

\pagestyle{emptyf}
\mbox{}

\end{document}